DOI: 10.1109/FAIML57028.2022.00048

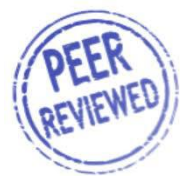

# Tawa Pukllay Proof: New Method for Solving Arithmetic Operations with The Inca Yupana Using Pattern Recognition and Parallelism

**Author:** Dhavit Prem (Universidad de Lima, Univ. Nacional Mayor de San Marcos)

**Coauthors:** Fernando Sotomayor (Universidad Nacional de Ingeniería) | Alvaro Saldívar (Universidad de Lima, Asociación Yupanki)| Rosario Guzman-Jimenez (Universidad de Lima)

**Abstract:**

Yupana is an Inca device used for arithmetic operations. This article describes a new arithmetical system: Tawa Pukllay (TP), where arithmetic operations do not require mental calculations: no carries, no borrows, no memorization of multiplication tables, nor trial and error procedures for divisions. Instead, user recognizes patterns and makes predefined movements to perform the four basic arithmetic operations very quickly; moreover, the result of the operation can be reached by multiple paths and in parallel, allowing each user to create his own strategies. This paper proves with mathematical rigor that TP produces correct numerical results.

# Tawa Pukllay Proof: New Method for Solving Arithmetic Operations with The Inca Yupana Using Pattern Recognition and Parallelism

Dhavit Prem*
Systems Engineering, University of Lima ULIMA
Master Philosophy-Epistemology
Universidad Nacional Mayor de San Marcos UNMSM
Asociacion Yupanki Lima, Peru
yupanki@yupanainka.com

Rosario Guzman-Jimenez
Faculty of Engineering
University of Lima ULIMA
Lima, Peru
rguzman@ulima.edu.pe

Fernando Sotomayor
Postgraduate School
Universidad Nacional de Ingenieria UNI
Lima, Peru
ffsotomayor@outlook.com

Alvaro Saldivar
Business Management, University of Lima ULIMA
Scientific and pedagogical research,
Asociacion Yupanki Lima, Peru
yachay@yupanainka.com

*Abstract*—Yupana is an Inca device used for arithmetic operations. This article describes a new arithmetical system: Tawa Pukllay (TP), where arithmetic operations do not require mental calculations: no carries, no borrows, no memorization of multiplication tables, nor trial and error procedures for divisions. Instead, user recognizes patterns and makes predefined movements to perform the four basic arithmetic operations very quickly; moreover, the result of the operation can be reached by multiple paths and in parallel, allowing each user to create his own strategies. This paper proves with mathematical rigor that TP produces correct numerical results.

## I. INTRODUCTION

The 16th century manuscript entitled "Nueva Corónica y Buen Gobierno" by Guaman Poma de Ayala contains a drawing that in its lower left corner shows a rectangle divided into 20 boxes organized in 5 rows and 4 columns, where the boxes of each row are marked from left to right with 5, 3, 2 and 1 points, as shown in Figure 1. The manuscript mentions that the rectangle represents the instrument that the Incas used in their accounting (1); but does not provide information on how this instrument was used.

Chronicler J. de Acosta wrote: *"But it seems a kinde of witchcraft, to see an other kinde of Quippos, which they make of* graines *of Mays(...) These Indians will take their graines and placefive of one side, three of another, and eight of another, and will change one graine of one side, and three of another. So as they finish a certaine account, without erring in any poynt: and they sooner submit themselves to reason by these Quippos, what every one ought to pay, than we cando with the penne."* (2)

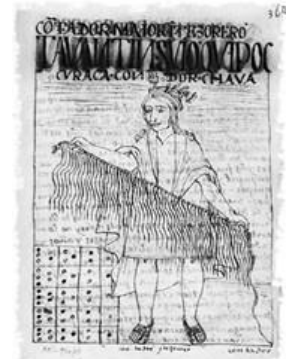

Figure 1. Cacique Condorchaua by Guaman Poma (3).

Since early 20th century, several hypotheses were formulated about how yupana (from Quechua "counting instrument")(4) was operated. These hypotheses differ mainly in four aspects (5): i) board orientation: one possible orientation is theone shown in Figure 1 (6; 7; 8; 9; 10; 11; 12; 13; 14), the other possible orientations are obtained by rotating the figure counterclockwise 90° or clock-wise 90° (15; 16; 17; 18); ii) interpretation of the dots: in one interpretation dots represent holes where to place the counting tokens (seeds or pebbles) (6; 16; 9; 11; 15; 13; 10; 18); in the other interpretation the number of dots indicate the numerical value assigned to each token deposited in the box (7; 17; 8; 12; 14); iii) numerical base: base 10 (6; 7; 8; 9; 10; 11; 12; 13; 16; 18; 15; 14; 19), binary base (20; 21), base 40 or multi-base (10); and iv) algorithms: ways of performing arithmetic operations (6; 7; 8; 9; 10; 11; 12; 13; 16; 18; 15; 14).

One of the researchers of the yupana, (8), formulated, among others, the following hypotheses: (i) the squares marked with the same number of dots are arranged in columns, as seen in Guaman Poma de Ayala's drawing; (ii) if a square is marked with *n* dots, that square is said to have weight *n*, which means that each token placed in this square can be replaced by *n* tokens placed in the square marked with 1 dot; iii) each row holds a digit of the decimal representation of a number; those rows from bottom to top are for units, tens, hundreds, thousands and tens of thousands, etc. respectively.

Tawa Pukllay (in Quechua tawa: four; pukllay: game) (4) is a method for performing arithmetic operations with the yupana developed by Dhavit Prem[1] and the Asociación Yupanki[2] (14). An important component of the method is a given set of pairs, each one made up of a pattern of squares with tokens and a corresponding token movement. The interesting and novel feature of the method is that the result of an arithmetic operation is achieved by repeatedly executing the following cycle:

i) visually detect on the yupana a pattern of squares with tokens, and ii) execute its corresponding token movement. The token movement corresponding to the detected pattern is mechanically executed, sequentially or in parallel, without having to do mental calculations such as carrying in addition, borrowing in subtraction, without having to remember multiplication tables, and without trial and error in division (5; 14; 22).

The method TP has three components: i) the representation of the numbers (taken from Pereyra's proposal), ii) the set of patterns of squares with tokens, along with their corresponding token movements (original proposal of TP), and iii) the instructions for performing the four arithmetic operations (original proposal of TP).

Some of the TP features are: i) if several patterns appear on the yupana, there is no specific order to execute their corresponding token movement; this allows to achieve the result by following different ways; ii) token movements can also be executed in parallel. (20); iii) TP looks like a game board such as chess, go, checkers, etc; this feature attracts students and encourages them to develop their own strategies (23; 24).

Patterns and their respective movements are held in a *Pattern-Move* Table (see Tables 1 and 2); each row of this table has a pattern and its associated movement.

TP has been presented at the Science and Technology Fair "Perú Conciencia" 2015-2017 (26; 25); Biblioteca Nacional del Perú (27) and Instituto Nacional de Cultura del Perú - Cusco (2016); USA Science & Engineering Festival in Washington DC (28), Latin American Congress of Mathematics RELME 31 - University of Lima (29), RELME 32 - University of Medellin, Colombia (30), VI International Congress of Ethnomathematics (31), EFPEM Congress at Universidad San Carlos de Guatemala (32), International Congress on Educational Innovation a the Instituto Tecnológico Superior de Monterrey, (Mexico 2019) (35).

Section 3.2 of this paper describes the instructions for arithmetic operations as high-level algorithms.

In Section 3.3 we publish for the first time the mathematical proof of TP correctness through theorems.

## II. DESCRIPTION OF THE YUPANA ACCORDING TO TP

Numbers are represented by placing small objects called **tokens** *(seeds, pebbles, etc.)* on the squares of the yupana. A square is identified by pair **(b,r)**, where **b** is the number of points, and **r** is the row where the square is.

### A. Representation of Numbers

*First rule of representation*

Say we represent in a yupana of $m$ rows a number that in base 10 is written as

$$\overline{d_{(m-1)}d_{(m-2)}\ldots d_2 d_1 d_0} \qquad (1)$$

Its digits from left to right are $d_{(m-1)}$, $d_{(m-2)}$, $d_2$, $d_1$ and $d_0$. The digit $d_0$ has to be represented in row 0 (the bottom one); digit $d_1$, in row 1; digit $d_2$, in row 2;...; and digit $d_{(m-1)}$, in row m-1 (the top row).(Figure 3).

According to this rule, row 0 is for the units ($10^0$); row 1, for the tens ($10^1$); row 2, for the hundreds ($10^2$), and so on.

This way of representing numbers determines that a yupana of $m$ rows can represent, besides 0, any number from 1 to $m$ decimal digits, that is, from 1 to $10^{m-1}$. Therefore, the yupana of Guaman Poma de Ayala, which has 5 rows, can represent from 0 to 99999.

To represent a decimal digit in a row, it is necessary to keep in mind that its squares are marked from left to right with 5, 3, 2 and 1 dots, and that the number of dots is called **weight of the square**.

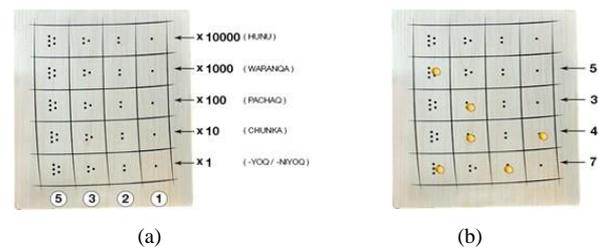

(a) (b)

Figure 2. (a)Yupana value structure (5) (b) Represents number 5347.

*Second rule of representation*

To represent a decimal digit d (from 1 to 9) in a row of the yupana, look for the smallest number of squares whose sum of weights is d, and place 1 token on each of these squares. Digit 0 is represented by leaving the row without tokens. Digits 1, 2, 3 and 5 are represented by a single token, which is placed in the square with the corresponding weight. Digit 4 is represented with one token in squares 1 and 3; digit 6, with one token in squares 5 and 1; digit 7, with one token in squares 5 and 2; digit 8, with 1 token in squares 5 and 3; digit 9, with 1 token in squares 5, 3 and 1 (see figure 2b).

*Third rule of representation*

For subtraction, minuend and subtrahend use different color of tokens. For division, dividend and divisor use different color of tokens too. Minuend and dividend tokens are called positive tokens, subtrahend and divisor tokens are called negative tokens.

### B. States and Values of the Yupana State of Yupana

Every possible configuration of tokens on the yupana is a state. A state is defined by the positive and negative tokens contained on each square.

*Simple state*

A state is the numeric value of the yupana on any position. A simple state is the optimal representation of that number, with the minimum quantity of tokens and the minimum quantity of squares.

*Weighted value of a square*

The weighted value of a square (b,r) is given by $b \times 10^r$.

*Effective value of a token*

The effective value of a positive token, is the weighted value of the square it occupies; the effective value of a negative token is the negative of the effective value of the square. That is, in a square (b,r) each positive token has effective value $b \times 10^r$; while each negative token has effective value $-b \times 10^r$.

*Effective value of a square*

The effective value of a square is the sum of the effective values of the tokens on it. In a given state, for each square(b,r) its number of positive tokens is $N_{b,r}^+$ and its number of negative tokens is $N_{b,r}^-$; therefore, the effective value of the square (b,r), which is denoted by $E_{(b,r)}$, is given by

$$E_{b,r} = \left(N_{b,r}^+ - N_{b,r}^-\right)\left(b \times 10^r\right) \qquad (2)$$

From the previous paragraph we see that an empty square ($N_{b,r}^+ = N_{b,r}^-$, r= 0) has effective value 0, and that for each positive token added to the square, the effective value of the square is increased by a value equal to the effective weight of the square, and that for each negative token added to the square, its effective value is decreased by the effective value of the square. For example, if a square of weight 3 in row 2 has 5 positive tokens added to it, its effective value is increased by $5(3 \cdot 10^2) = 5 \cdot 300 = 1500$

*Value of the yupana*

(**value represented by the yupana**). It is defined as the sum of the effective values of all the squares of the yupana. Recalling that the effective value of each square *(b,r)* is denoted by $E_{(b,r)}$, the value of the yupana is given by $\hat{V} = \sum_{(b,r) \in C} E_{b,r}$, where C is the set of squares of the yupana.

The value of the yupana can also be seen as the sum of the effective values of all the *tokens* that are in the yupana. This is because $E_{(b,r)}$ is the sum of effective values of the tokens in the square *(b,r)*, so that

$$\hat{V} = \sum_{(b,r) \in C} E_{b,r} = \sum_{(b,r) \in C} \textit{effective value tokens } (b,r) \qquad (3)$$

This is the main formula of the yupana according to TP. All its mathematical properties are derived from it.

When the yupana represents a value *Y*, it is said that *the value of the yupana is Y*, or that *the yupana has the value Y*, or that *the yupana contains the value Y*. When tokens are deposited so that the yupana represents a desired value Y, it is said that *the value Y is deposited*, or that *the value Y is stored*, or that *the value Y is loaded*.

## III. ARITHMETIC OPERATIONS WITH TAWA PUKLLAY (TP)

Arithmetic operations with TP are performed using the **pattern-movement table** (Tables 1,2). There are two types of actions: i) **load operands** into the yupana and ii) perform **simplification**. Loading an operand consists of placing tokens on the yupana to represent an operand according to the rules of representation of a number. Simplification consists of **visually detect** on the yupana a **pattern** from the pattern-movement table and **execute** the movement associated with the **detected pattern**, until answer is found in simple state and there are no operands left to load. The movements in the table are of three types:

- **Reducing (basic) movements**, which always decrease the number of tokens in the yupana or in a square of the yupana (without altering the value); with these movements, the value of the yupana is expressed with the least number of tokens so that it can be recognized at a glance.
- **Expansion movements**, which are the inverse of the reducing movements and increase the number of tokens in the yupana (without altering their value); these movements are necessary for subtraction and division (they can also be used for addition and multiplication). With them, when there are positive and negative tokens, each negative token is paired with a positive token, and then each pair of a positive token with a negative token is eliminated.
- **Composite movements**, which are equivalent to the execution of two or more movements of any type. They are used to perform a sequence of movements with a single movement.

*Simplification of the yupana*

In every state the yupana represents a value; however, that value is not always recognizable at first glance. To be recognized at a glance, the value must be expressed with as few tokens as possible, i.e., the yupana must be in a simple state. This is done by executing movements from the pattern-movement table. This process of carrying the yupana to a simple state without changing its value is called simplification of the yupana.

An example of simplification is the following. In the cases of addition, subtraction and multiplication, just by loading the operands on the yupana, the value of the yupana is the result of the operation; but this value may not be recognizable at simple glance, and to make it recognizable one has to execute movements of the pattern-movement table that, without changing its value, bring the yupana to a simple state. This is all that has to be done to obtain the result of the addition, subtraction and multiplication.

*Non-determinism and parallelism*

It is possible that in a given state, several patterns of the pattern-movement table may appear in various parts of the yupana. When this occurs, any of the present patterns can be chosen and its corresponding movement can be executed; the movements could even be executed in parallel.

A. *Specific Instructions for Arithmetic Operations with the Yupana*

1) *Instructions to perform addition ("Yapay"):* Two or more summands can be added together (Figure 3):

a) *Load the addends:* load all the addends one on top of the other; that is, each addend is loaded on top of the ones already loaded. At the end of the loading of the addends, the value that the yupana has is the result of the addition; however, this value may not be recognizable at simple glance (simple state).

*b) Simplify the result of the addition:* if the result of the addition is not recognizable at simple glance, execute the cycle detect pattern-execute move until the number contained in the yupana is in simple state.

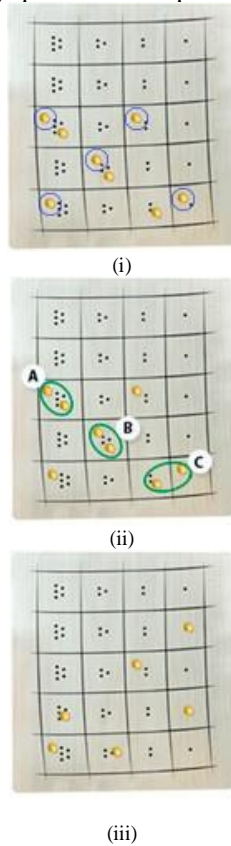

Figure 3.  736+532 (i) First addend in blue circles. Second addend with the other tokens (ii) Patterns A, B, C are detected and their movements are executed (iii) Step ii repeats until each row represents a digit in simple state. Answer=1268.

*2) Instructions to perform subtraction ("T'aqay"):* Subtraction can be performed with several minuends and subtrahends loaded in the yupana in any order (Figure 4):

*a) Load the minuends and subtrahends:* load the minuends and subtrahends in any order (they can be intermixed), the minuends with tokens of one color and the subtrahends with tokens of a different color. The minuend tokens are called positive tokens, and the subtrahend tokens are called negative tokens. When the operands are completely loaded, the value of the yupana is the result of the subtraction, which may not be recognizable at first glance.

*b) Simplify the result of subtraction:* In the case of subtraction, the first goal of simplification is to have the yupana with tokens of a single color. Proceed as follows:
- Remove all pairs of tokens of different color on the same square.
- If the yupana still has tokens of both colors, execute the *cycle detect pattern - execute movement* to bring positive and negative tokens to same square, then go back to step a.

- Once all the tokens on the yupana are of a single color (positive or negative), continue with simplification until the value is in simple state.

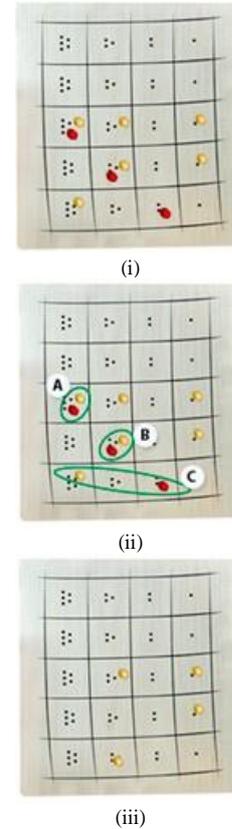

Figure 4.  945 - 532 (i) Minuend with yellow, subtrahend with red (ii) We look forward to pair yellow and red tokens on same squares, so patterns A, B, C are detected and movements are performed (iii) We use expansion movements repeating step ii until there are only tokens of one color representing the answer in simple state. Answer=413.

*3) Instructions to perform the multiplication ("Miray"):* Multiplication is interpreted as an addition in which the multiplicand intervenes adding itself the number of times indicated by the multiplier. (Figure 5).

**Replication of a token**

Replication of *n* times a token consists of replacing that token with *n* tokens in its own square. Replication is necessary in multiplication where the multiplier is *n*. If *n* is large, replication would require a large number of tokens. Therefore, to achieve the same effect that replication produces, but using fewer tokens, the following **abbreviated replication of a token** procedure is used, which is presented below.

**Procedure for replicating a token of the multiplicand.**

Suppose that the multiplier is *n*, which in decimal notation is written as $\overline{d_k\ d_{(k-1)}\ d_{(k-2)}...d_2 d_1 d_0}$

- Place $d_k$ tokens in the same column as the token being replicated, but *k* rows higher. Place $d_{(k-1)}$ tokens in the same column as the token being replicated, but *k-1* rows higher. Continue in this manner with the digits $d_{k-2}, ..., d_2, d_1, d_0$.
- Remove from the yupana the token of the multiplicand just replicated.

**Abbreviated replication of a token**

The abbreviated replication of *n* times a token produces the same change in the value of the yupana as the *replication of n times a token* produces, but uses fewer tokens and movements.

**Procedure for multiplication**

- **Load the multiplicand.**
- **Replicate the multiplicand the number of times indicated by the multiplier.** Let *n* be the multiplier. Replicate *n* times each token of the multiplicand using the procedure **abbreviated replication of a token** presented in previous paragraph. When all the tokens of the multiplicand have been replicated, the value of the yupana is the result of the multiplication; however, this value may not be recognizable at first glance.
- **Simplify the result of the multiplication:** if the result of the multiplication is not recognizable at first glance, execute the cycle *detect pattern-execute movement* until the number contained in the yupana is recognizable at a glance (simple state).

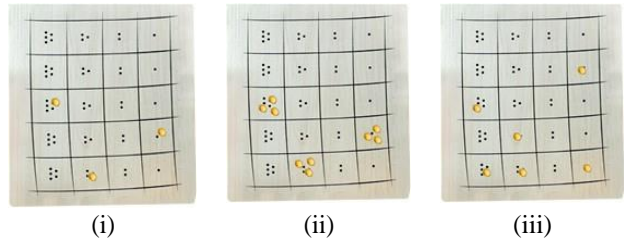

(i)      (ii)      (iii)

Figure 5: Example: 513 x 3 (i) Represent the multiplicand 513 (ii) Replicate the multiplicand as many times as pointed by the multiplier 3 (iii) Recognize patterns and execute movements until we get the answer in simple state. Answer = 1539

*4) Instructions to perform the division ("Rakiy"):* The quotient of the division is interpreted as the number of times the dividend contains the divisor. To find this number, the divisor is subtracted from the dividend as many times as necessary until the remainder of the dividend is a value less than the divisor. Figure 6 and Figure 7 show an example of this process.

**Repeated subtraction.** Conceptually, it consists in subtracting the divisor from the dividend; then from the result subtract again the divisor; and so on, repeatedly until the result of the subtraction is less than the divisor. The result of the last subtraction is the remainder of the division, and the number of times the subtraction was performed is the quotient. However, if the ratio of the dividend to the divisor is large, the above process may require a large number of subtractions. So, instead of using repeated subtraction it has to be used **fast repeated subtraction** given next.

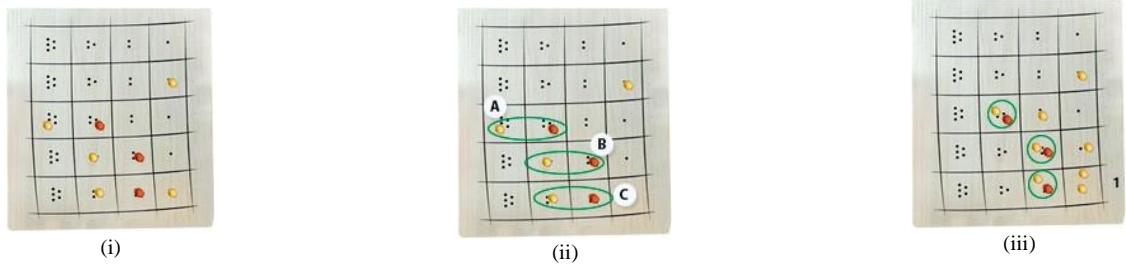

(i)      (ii)      (iii)

Figure 6: Example: 1534/322 (i) Yellow tokens for dividend, red tokens for divisor (ii) Applying movements from the pattern-movements table, we pair yellow and red tokens on the same squares (iii) When each red token is paired with a yellow one, the yellow ones are taken away and subtraction counter increases by 1.

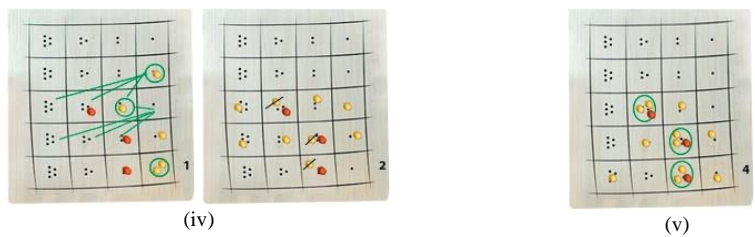

(iv)      (v)

Figure 7: (iv) Repeat steps ii and iii until there is no yellow token or until the dividend (yellow tokens) is less than the divisor (red tokens) (v) The number of times the yellow tokens have been removed represents the quotient of the division: 4. The remaining yellow tokens represent the remainder: 246.

**Fast repeated subtraction** ($S_a$). It is a sequence of movements that produces the same effect as repeated subtraction, but with fewer movements. The steps are as follows:

- Initialize the subtraction counter $q : q \leftarrow 0$; this counter will indicate how many times the divisor has been subtracted from the dividend, so its final value will be the quotient.

- Displace all the divisor tokens in their own columns by the same number of rows until the displaced divisor is the closest value to the current dividend but not greater than it; initialize the displace indicator

$k : k \leftarrow$ *quantity of displaced rows*

- Subtract the displaced divisor from current dividend: Match each token of the *displaced divisor* (negative token) with a token of the current dividend (positive token) by making movements from the *pattern-movement* table.
  From each of the pairs remove only the token of the current dividend, the token of the displaced divisor remains in the yupana to perform the next subtraction. Subtraction counter increments in $10^k$: $q \leftarrow q + 10^k$.
- If the current dividend (obtained in step c) is less than the *displaced divisor*, move the *displaced divisor* down until it is less or equal than the *current dividend*; for each row displaced down, the *displace indicator* has to be updated: $k \rightarrow k-1$; Then we go back to step *c*).

If the *displaced divisor* can no longer be moved down, it means that the end of the *repeated subtraction* has been reached. The number of times the divisor has been subtracted is in the *subtraction counter q*, which is the *quotient of the division*, and *the temporary remainder* is the *remainder* of the division, which may not be in simple state.

**Procedure for division.**
- **Load operands:** load the dividend with positive tokens, and the divisor with negative tokens.
- **Perform the Fast Repeated Subtraction:** conceptually it consists of the following: subtract the divisor from the dividend; the divisor is subtracted again from the result of the subtraction and so on repeatedly until the result of the subtraction is less than the divisor. The result of the last subtraction is the remainder of the division, and the number of times the subtraction was made is the quotient. The *temporary remainder* is the *remainder* of the division, which may not be in simple state.
- **Simplify the remainder:** if the remainder is not recognizable at first glance, execute the cycle *detect pattern-execute move* until the number contained in the yupana is recognizable at a glance (simple state).

*B. Mathematical Validity of TP*

*Theorem 1:* **Transfer:** If the effective value of a square has a variation $\delta$ (positive or negative), the value of the yupana also gets a variation $\delta$.
*Proof:*
Suppose that square $(\sigma, \rho)$ has effective value $E_{\sigma,\rho}$ and that the yupana has the value $Y$. By definition,

$$Y = \sum_{(b,r) \in C} E_{(b,r)}$$

where $E_{(b,r)}$ is the effective value of square $(b,r)$ and C is the set of squares of the yupana. Extracting from the summation the effective value of the square $(\sigma, \rho)$ we have

$$Y = E_{(\sigma,\rho)} + \sum_{(b,r) \in C - \{(\sigma,\rho)\}} E_{b,r} \quad (4)$$

Suppose now that $E_{\sigma,\rho}$ changes to $(E_{\sigma,\rho} + \delta)$, the new value of the yupana will be

$$Y' = \left( (E_{\sigma,\rho} + \delta) + \sum_{(b,r) \in C - \{(\sigma,\rho)\}} E_{b,r} \right) = \left( E_{\sigma,\rho} + \sum_{(b,r) \in C - \{(\sigma,\rho)\}} E_{b,r} \right) + \delta = Y + \delta \quad (5)$$

*Theorem 2:* **Invariance:** The value of a yupana is invariant under movements of pattern-movement table
*Proof:* Appendix contains a *pattern-move* table with movements containing the corresponding token changes (withdrawal or deposit) made on the involved squares, together with the variation that these changes produce in the effective value of the square. The algebraic sum of the variations of the effective value of the involved squares is zero. Therefore, due to Thm. 1, the variation of the value of the yupana is 0.

*Theorem 3:* **Superposition:** When two or more numbers are superimposed on the yupana, the value of the yupana is the sum of the superimposed numbers.
*Proof:* Suppose we have two yupanas, A and B, with values $Y_A$ and $Y_B$. Suppose also that each square $(b,r)$ of A has effective value $E^{(A)}_{(b,r)}$ and the corresponding square of B has $E^{(B)}_{(b,r)}$. Then we have:

$$Y_A = \sum_{(b,r) \in C} E^{(A)}_{(b,r)} \text{ and } Y_B = \sum_{(b,r) \in C} E^{(B)}_{(b,r)} \quad (6)$$

where C is the set of squares of both A and B. Now, if all tokens on square $(b,r)$ of yupana B are moved to square $(b,r)$ of yupana A, the effective value of square $(b,r)$ of yupana A will be $E^{(A)}_{(b,r)} + E^{(B)}_{(b,r)}$, and the new value of yupana A will be

$$Y'_A = \sum_{(b,r) \in C} \left( E^{(A)}_{(b,r)} + E^{(B)}_{(b,r)} \right) = \left( \sum_{(b,r) \in C} E^{(A)}_{(b,r)} \right) + \left( \sum_{(b,r) \in C} E^{(B)}_{(b,r)} \right) = Y_A + Y_B \quad (7)$$

*Theorem 4:* If every token of a yupana which is in a simple state is replaced on its own square by $k$ tokens, the value of the yupana is multiplied by $k$.
*Proof:* Suppose that the yupana is in a simple state and has a value Y. The action of replacing each token by $k$ tokens on its own square is equivalent to superimposing $k$ times the same number Y, but by Thm. 3 the result of the superposition is the sum of the superimposed numbers, so the number Y is added $k$ times: $kY$.

*Theorem 5:* Replication of a token $n$ times and abbreviated replication of the same token $n$ times produce the same change in the yupana value.
*Proof:*
Suppose that n in decimal notation is written as $\overline{(d_k d_{k-1} d_{k-2} \ldots d_2 d_1 d_0)}$, then

$$n = d_k 10^k + d_{k-1} 10^{k-1} + d_{k-2} 10^{k-2} + \ldots + d_2 10^2 + d_1 10 + d_0 = \sum_{i=k,k-1,\ldots,1,0} d_i 10^i \quad (8)$$

The step 1 of the *abbreviated replication* of $n$ times a token of the square $(b,r)$ consists in for each $i = k, k-1,\ldots, 2, 1, 0$ $d_i$ tokens are deposited in the square of column $b$ row $(r+i)$, whereby the effective value of square $(b,r+i)$ varies by $d_i(b \times 10^{r+i})$. Because of Thm. 1 this variation is transferred to the value of the yupana. Therefore, the total variation in the value of the yupana is

$$\sum_{i=k,k-1,\ldots,1,0} d_i (b \cdot 10^{r+1}) = \sum_{i=k,k-1,\ldots,1,0} d_i 10^i (b \cdot 10^r) = (b \cdot 10^r) \sum_{i=k,k-1,\ldots,1,0} d_i 10^i$$
$$(b \cdot 10^r) \cdot n \quad (9)$$

Step 2 of the *abbreviated replication* consists of removing the replicated token from its square $(b,r)$, which changes the effective value of the square by

$$-(b \cdot 10^r) \quad (10)$$

Therefore, considering (Eq.10) and (Eq.9) the variation in the yupana value produced by the *abbreviated replication* of $n$ times the token in the square $(b,r)$ is

$$(b \cdot 10^r) \cdot n - (b \cdot 10^r) = (n-1) \cdot (b \cdot 10^r) \quad (11)$$

On the other hand, replicating $n$ times the token in square $(b,r)$ is equivalent to superimposing $n-1$ tokens on the token in that square, and therefore the variation in the effective value of that square is $(n-1) \cdot (b \cdot 10^r)$ (Eq.10).

Finally, by (Eq.10) and (Eq.11) it is shown that the *abbreviated replication* of $n$ times a token and the *replication* of $n$ times produce the same variation in the yupana value.

*Theorem 6*: **Validity of addition:** The result of the addition performed according to TP is correct.

*Proof:* When superimposing the summands, by Thm. 3 *(superposition)*, the value of the yupana is the sum $S$ of the summands. In order to have $S$ in simple state, simplification is performed which, by Thm. 2 *(invariance)*, leaves the yupana with the value $S$. Therefore, the result of the sum performed according to TP is correct.

*Theorem 7*: **Validity of subtraction:** The result of the subtraction performed according to TP is correct.

*Proof:*

When superimposing the minuends and subtrahends, by the Thm. 3 *(superposition)*, the value of the yupana is the difference between the sum of all the minuends minus the sum of all the subtrahends, let's call $D$ such a value. In order to have $D$ in simple state, simplification is performed which, by Thm. 2 *(invariance)*, leaves the yupana with the value $D$. Therefore, the result of subtraction performed according to TP is correct.

*Theorem 8*: **Validity of multiplication:** The result of multiplication performed according to TP is correct.

*Proof:*

Let $m$ be the multiplicand and $n$ the multiplier. The replication of $m$, which consists of replacing each of its tokens by $n$ tokens on its same square, is equivalent to superimposing the multiplicand $m$ on itself $n$ times. By Thm. 3 *(superposition)*, superimposing the multiplicand $m$ on itself $n$ times is equivalent to the sum where the multiplicand $m$ appears $n$ times, i.e., $(m+m+\ldots+m)$ "$n$ times" but this is $m \times n$. Then, after the replication of the multiplicand, the yupana has the result of the multiplication $m \times n$. If the value $m \times n$ is not in simple state then simplification is performed, and by Thm. 2 *(invariance)*, leaves the yupana with value $m \times n$. Therefore, the result of the multiplication performed according to TP is correct. The method does not use the replication directly, instead it uses the abbreviated replication, but by Thm. 5 they produce the same result.

*Theorem 9*: **Validity of division:** The result of the division performed according by TP is correct.

*Proof:*

The result of the division is made by subtracting the divisor from the dividend repeatedly until a remainder smaller than the divisor is obtained. The number of times the subtraction is performed is the quotient and the remainder of the dividend is the remainder. Due to Thm. 7 *(validity of subtraction)*, the result of the division performed according to TP is correct.

IV. DISCUSSION

The purpose of this paper was to demonstrate the mathematical validity of TP. The demonstration is based on 4 facts:

i) that in the case of addition and subtraction, when loading the operands, the value of the yupana is equal to the value of the result, ii) that the movements of the *pattern-move* tables do not alter the value of the yupana, iii) that multiplication is interpreted as a repeated addition, iv) that division is interpreted as a repeated subtraction of the divisor from the dividend.

TP, unlike traditional methods, allows the representation of numbers in a simple way that facilitates their reading. This is due to the structure of dots and squares of the yupana that allow the recognition of digits and numbers in an immediate way (subitization).

Arithmetic operations with many numbers often demand stressful effort. In contrast, TP presents a playful procedure as a strategy-oriented challenge, which reduces response time compared to traditional methods. Furthermore, with TP it is possible to execute one operation by two or more students on the same board in parallel. For this, we would divide the board in zones and assign them to each student. This allows to develop mathematical strategies in parallel and to introduce intuitively collaborative work. (35).

TP represents zero as number on the yupana by leaving the board without tokens, meanwhile zero as representation of digit on any decimal place value can be represented by leaving a row without tokens. This type of representation of zero by leaving spaces in blank can be also seen in the *khipus*, the inca knots to store data. (36; 16; 11; 19).

Carry-overs (addition and multiplication), borrowings (subtraction), multiplication tables memorization (multiplication and division) and trial and error (division) are frequent difficulties when using traditional methods. TP eliminates all those difficulties by offering a simple pattern-movement procedure.

With TP subtractions can be calculated with several minuends and subtrahends in the same operation. In contrast, traditional methods can only perform subtractions with just one minuend and one subtrahend.

Using TP for divisions, students can easily see the meaning of division; in contrast, with traditional methods, students apply trial and error, often without understanding what they are doing.

By seeing that addition and multiplication operands can be loaded in different order, and that a quantity can be represented in many ways, students intuitively assimilate the *commutativity* and *associativity* properties.

TP is non-deterministic, this means that when several patterns from the *pattern-move* table are detected on the yupana, we can choose any of them, and execute its associated movement. This TP feature catches the attention of students, and offers them a chance to design their own pattern-movements and optimal strategies.

A playful way to introduce TP is the ancestral andean practice called *Atipanakuy*, (in quechua *"mutual empowerment"*) which is a competition of two or more students performing arithmetic challenges in the shortest time; this encourages them to create better strategies (24; 33; 35).

## V. CONCLUSIONS

This paper has demonstrated with mathematical rigor that TP can solve the four basic arithmetic operations with integers by pattern recognition, executing movements in parallel and without use of traditional indoarabic calculations. It has also been shown that TP features advantages to overcome frequent difficulties in the arithmetic teaching and learning process *(carry-overs, borrows, multiplication tables, division trial and error procedure)*. Next papers will demonstrate that: a) TP can solve the four basic arithmetical operations also with decimal fractions and other arithmetical operations such as power, square root, percentage, logarithm, etc. with integers b) TP can be used with multi-value tokens, and c) TP logic can be used in software and hardware development.

## VI. APPENDIX

- **Reducing moves:** always decrease the number of tokens without altering the value of the yupana.
- **Expansion moves:** always increase the number of tokens in the yupana without altering their value.
- **Composite moves:** they are the combination of two or more movements of any type.

Note 1: In this table, [5],[3],[2] and [1] are the names of the squares of the yupana.

Note 2: In the formulas of this table, terms preceded by a "+" sign represent the variation when tokens are placed on the board, terms preceded by a "-" sign represent the variation when tokens are removed from the board. As the formulas show, variations of the value of the yupana are always zero.

TABLE I. PATTERNS-MOVEMENTS TABLE: REDUCING MOVEMENTS (A.K.A. BASIC MOVEMENTS)

| Pattern → Movement | Movement description and variation verification formula |
|---|---|
| (image) | **Iskay "Short opening":** [2] has more than one token<br>**Movement:** a) If the number of tokens on [2] is even, move one half to [1]; and the other half to [3]. b) If the number of tokens on [2] is odd, leave one token on [2]; apply procedure a)<br>**Variation:** $(-2k - 2k + 3k + k) \cdot 10^r = 0$ |
| (image) | **Kimsa "Long opening":** [3] has more than one token<br>**Movement:** a) If number of tokens is even, move one half to [1]; and the other half to [5]. b) If number of tokens is odd, leave one token on [3]; apply procedure a)<br>**Variation:** $(-3k - 3k + 5k + k) \cdot 10^r = 0$ |
| (image) | **Pisqa "Big L":** [5] at row $r$ has more than one token<br>**Movement:** a) If number of tokens is even, move half of them to [1] of row $r+1$; remove the other half from the board. b) If number of tokens is odd, leave 1 token on [5]; apply procedure a) with remaining tokens. **Variation:** $(-5k - 5k) \cdot 10^r + k \cdot 10^{r+1} = 0$ |
| (image) | **Kikin "2 on [1] is 1 on [2]":** When [1] has two tokens<br>**Movement:** Replace two tokens on [1] by one token on [2]<br>**Variation:** $(-k - k + 2k) \cdot 10^r = 0$ |
| (image) | **Kikin "3 on [1] is 1 on [3]":** When [1] has three tokens<br>**Movement:** Replace three tokens on [1] by one token on [3]<br>**Variation:** $(-k - k - k + 3k) \cdot 10^r = 0$ |
| (image) | **Kikin "5 on [1] is 1 on [5]":** When [1] has 5 tokens<br>**Movement:** Replace five tokens on [1] by one token on [5]<br>**Variation:** $(-k - k - k - k - k + 5k) \cdot 10^r = 0$ |
| (image) | **Pichana [1][2] "Sweep [1][2]":** When [1] and [2] have both one token each<br>**Movement:** Replace both tokens by one token on [3]<br>**Variation:** $(-k - 2k + 3k) \cdot 10^r = 0$ |
| (image) | **Pichana [2][3] "Sweep [2][3]":** When [2] and [3] have both one token each<br>**Movement:** Replace both tokens by one token on [5]<br>**Variation:** $(-2k - 3k + 5k) \cdot 10^r = 0$ |

TABLE II.  PATTERNS-MOVEMENTS TABLE: EXPANSION & COMPOSITE MOVEMENTS

| Pattern → Movement | Movement description and variation verification formula |
|---|---|
| → | **Expansion of 5:** One token on [5]<br>**Movement:** Replace 1 token on [5] by 1 token on [3] and 1 token on [2].<br>**Variation:** $(-5k + 3k + 2k) \cdot 10^r = 0$ |
| → | **Expansion of 3:** One token on [3]<br>**Movement:** Replace 1 token on [3] by 1 token on [2] and 1 token on [1]<br>**Variation:** $(-3k + 2k + k) \cdot 10^r = 0$ |
| → | **Expansion of 2:** One token on [2]<br>**Movement:** Replace 1 token on [2] by 2 tokens on [1]<br>**Variation:** $(-2k + k + k) \cdot 10^r = 0$ |
| → | **Inverse Pisqa:** One token on [1] of row $r+1$<br>**Movement:** Replace the token on [1] of row $r+1$ by 2 tokens on [5] of row $r$<br>**Variation:** $-k \cdot 10^{r+1} + (5k + 5k) \cdot 10^r = 0$ |
| → | **Inverse Hatun Pichana:** One token on [1] of row $r+1$<br>**Movement:** Replace that token by 1 token on [5], 1 token on [3] and 1 token on [2] of row $r$<br>**Variation:** $-k \cdot 10^{r+1} + (5k + 3k + 2k) \cdot 10^r = 0$ |
| → | **Inverse Sonqo ("inverse heart"):** One token on [1] of row $r+1$<br>**Movement:** Replace that token by 2 tokens on [3] and 2 tokens on [2] of row $r$.<br>**Variation:** $-k \cdot 10^{r+1} + (3k + 3k + 2k + 2k) \cdot 10^r = 0$ |
| → | **Inverse Huq Iskay Kimsa ("inverse 1-2-3"):** One token on [1] of row $r+1$<br>**Movement:** Replace that token by 1 token on [3], 2 tokens on [2] and 3 tokens on [1] of row $r$.<br>**Variation:** $-k \cdot 10^{r+1} + (3k + 2k + 2k + k + k + k) \cdot 10^r = 0$ |
| → | **Chunka "10":** *(Also Pachaq '100', Waranqa '1000', Hunu '10000',...)*<br>Any square **[s]** has $10^n$ tokens at row $r$ ($\forall n \in N; n \geq 1$)<br>**Movement:** Replace $10^n$ tokens on [s] of row $r$ by 1 token on [s], of row $r+n$.<br>**Variation:** $-s \cdot 10^n \cdot (10^r) + s \cdot (10^{r+n}) = 0$ |
| → | **Sonqo ("Heart"):** Two tokens on [2] and two tokens on [3] of row $r$<br>**Movement:** Replace those four tokens by 1 token on [1] of row $r+1$.<br>**Variation:** $(-3k - 3k - 2k - 2k) \cdot 10^r + k \cdot 10^{r+1} = 0$ |
| → | **Hatun Pichana ("Big sweep"):** One token on [5], one token on [3] and one token on [2] of row r<br>**Movement:** Replace those 3 tokens by 1 token on [1] of row $r+1$.<br>**Variation:** $(-5k - 3k - 2k) \cdot 10^r + k \cdot 10^{r+1} = 0$ |
| → | **Paña Chaska ("Star at right side"):** Two tokens on [3], one token on [2], two tokens on [1] of row $r$<br>**Movement:** Replace all those 5 tokens by 1 token in square [1] of row $r+1$<br>**Variation:** $(-3k - 3k - 2k - k - k) \cdot 10^r + k \cdot 10^{r+1} = 0$ |
| → | **Huq-Iskay-Kimsa ("1-2-3"):** One token on [3], two tokens on [2], three tokens on [1] of row $r$<br>**Movement:** Replace those six tokens by 1 token on [1] of row $r+1$.<br>**Variation:** $(-3k - 2k - 2k - k - k - k) \cdot 10^r + k \cdot 10^{r+1} = 0$ |
| → | **K'usillu ("Monkey"):** Three tokens on [3] and one token on [1] of row $r$<br>**Movement:** Replace those four tokens by 1 token in [1] of row $r+1$.<br>**Variation:** $(-3k - 3k - 3k - k) \cdot 10^r + k \cdot 10^{r+1} = 0$ |
| → | **Chinkay ("Dissapear"):** $k$ positive tokens and $k$ negative tokens on [s]<br>**Movement:** Take away k pairs of both positive and negative tokens of [s].<br>**Variation:** $-k(s \cdot 10^r) + k(s \cdot 10^r) = 0$ |

NOTES

[1]This article will be part of the Master's Thesis at the Faculty of Philosophy with mention in Epistemology at Universidad Nacional Mayor de SanMarcos (Lima, Perú) presented by Dhavit Prem (a.k.a. Carlos Saldívar Olazo).

[2]For more information about this research, you can contact the website www.yupanainka.com or send an email to author or coauthors.